\documentclass[leqno,11pt,a4paper]{amsart}
\usepackage{graphicx}
\usepackage[usenames,dvipsnames]{color}
\usepackage{hyperref}
\usepackage{subfigure}
\usepackage{array,ragged2e}
\usepackage{mathrsfs}
\usepackage{tabularx}
\usepackage{algpseudocode}
\usepackage{mleftright}
\usepackage{amssymb}
\usepackage{booktabs}
\usepackage{colortbl}
\usepackage{array}
\newtheorem{thm}{Theorem}[section]

\theoremstyle{definition}
\newtheorem{example}[thm]{Example}
\newtheorem{remark}[thm]{Remark}

\newtheorem{definition}[thm]{Definition}

\numberwithin{equation}{section}
\long\def\blankfootnotetext#1{\begingroup\def\thefootnote{\fnsymbol{footnote}}\footnotetext{#1}\endgroup}
\newcommand{\pro}[2]{\langle{#1},{#2}\rangle}
\newcommand{\Hom}{\mathrm{Hom}}

\newcommand{\Z}{\mathbb{Z}}
\newcommand{\Q}{\mathbb{Q}}

\newcommand{\C}{\mathbb{C}}
\newcommand{\Proj}{\mathbb{P}}
\newcommand{\V}[1]{\mathcal{V}\mleft({#1}\mright)}

\newcommand{\GL}{\mathrm{GL}}
\newcommand{\Aut}[1]{\mathrm{Aut}\mleft({#1}\mright)}
\newcommand{\AffAut}[1]{\mathrm{AffAut}\mleft({#1}\mright)}
\newcommand{\G}[1]{\mathscr{G}\mleft({#1}\mright)}
\newcommand{\magma}{{\sc Magma}}
\newcommand{\sage}{{\sc Sage}}
\newcommand{\palp}{\textsc{Palp}}
\newcommand{\idx}[1]{\abs{{#1}:\Lambda}}
\newcommand{\Vol}[1]{\mathrm{Vol}\mleft({#1}\mright)}
\newcommand{\linspan}[1]{\mathrm{span}\mleft({#1}\mright)}
\newcommand{\aff}[1]{\mathrm{aff}\mleft({#1}\mright)}
\newcommand{\cone}[1]{\mathrm{cone}\mleft({#1}\mright)}
\newcommand{\conv}[1]{\mathrm{conv}\mleft({#1}\mright)}
\newcommand{\sconv}[1]{\mathrm{conv}\mleft\{{#1}\mright\}}
\newcommand{\Newt}[1]{\mathrm{Newt}\mleft({#1}\mright)}
\newcommand{\NF}[1]{\mathrm{NF}\mleft({#1}\mright)}
\newcommand{\AffNF}[1]{\mathrm{AffNF}\mleft({#1}\mright)}
\newcommand{\coeffs}[1]{\mathrm{coeffs}\mleft({#1}\mright)}
\newcommand{\abs}[1]{\left\vert{#1}\right\vert}
\newcommand{\mmax}[1]{#1^\text{max}}
\newcommand{\PM}{P\!M}
\newcommand{\PMmax}{\mmax{\PM}}
\newcommand{\IndexOfMax}[1]{\mathrm{IndexOfMax}\mleft\{{#1}\mright\}}
\newcommand{\Continue}{\State\textbf{continue}}
\newcommand{\Break}{\State\textbf{break}}
\newcommand{\ForTo}{\textbf{to}\ }
\newcommand{\By}{\textbf{by}\ }
\newcommand{\SuchThat}{\textbf{such that}\ }
\graphicspath{{images/}}
\begin{document}
\author[R.~Grinis]{Roland Grinis}
\address{Trinity College\\University of Cambridge\\Cambridge, CB$2$\ $1$TQ\\UK}
\email{roland.grinis09@imperial.ac.uk}
\author[A.~M.~Kasprzyk]{Alexander M.~Kasprzyk}
\address{Department of Mathematics\\Imperial College London\\London, SW$7$\ $2$AZ\\UK}
\email{a.m.kasprzyk@imperial.ac.uk}
\blankfootnotetext{\emph{Keywords}: Convex polytope, computation, isomorphism, automorphism, normal form, {\palp}.}
\blankfootnotetext{2010 \emph{Mathematics Subject Classification}: 52B20 (Primary); 52B55, 52C07 (Secondary).}
\title{Normal forms of convex lattice polytopes}
\begin{abstract}
We describe an algorithm for determining whether two convex polytopes $P$ and $Q$, embedded in a lattice, are isomorphic with respect to a lattice automorphism. We extend this to a method for determining if $P$ and $Q$ are equivalent, i.e.~whether there exists an affine lattice automorphism that sends $P$ to $Q$. Methods for calculating the automorphism group and affine automorphism group of $P$ are also described.

An alternative strategy is to determine a normal form such that $P$ and $Q$ are isomorphic if and only if their normal forms are equal. This is the approach adopted by Kreuzer and Skarke in their {\palp} software. We describe the Kreuzer--Skarke method in detail, and give an improved algorithm when $P$ has many symmetries. Numerous examples, plus two appendices containing detailed pseudo-code, should help with any future reimplementations of these techniques. We conclude by explaining how to define and calculate the normal form of a Laurent polynomial.
\end{abstract}
\maketitle
\section{Introduction}
Determining whether two convex polytopes $P$ and $Q$, embedded in a lattice $\Lambda$, are isomorphic with respect to a lattice automorphism is a fundamental computational problem. For example, in toric geometry lattice polytopes form one of the key constructions of projective toric varieties, and any classification must somehow address the issue of whether there exists an automorphism of the underlying lattice sending $P$ to $Q$. In general, any isomorphism problem can be solved in one of two ways: on a case-by-case basis by constructing an explicit isomorphism between the two objects, or by determining a normal form for each isomorphism class.

The first approach -- dynamically constructing a lattice-preserving isomorphism in $\GL_n(\Z)$ between the two polytopes -- is discussed in~\S\ref{sec:iso_via_face_graph}. We describe one possible way to determine isomorphism of polytopes via the labelled face graph $\G{P}$ (see~\S\ref{subsec:labelled_face_graph}). This has the advantage that it works equally well for rational polytopes and for polytopes of non-zero codimension. By reducing the problem to a graph isomorphism question, well-developed tools such as Brendan McKay's~\textsc{Nauty} software~\cite{McKay81,Nauty} can then be applied.

Because our approach to isomorphism testing works equally well for rational polytopes, we are able to answer when two polytopes are equivalent, i.e.~when there exists an isomorphism $B\in\GL_n(\Z)$ and lattice translation $c\in\Lambda$ such that $PB+c=Q$. This is discussed in~\S\ref{subsec:equivalence}.  We can also calculate the automorphism group $\Aut{P}\leq\GL_n(\Z)$ of $P$: this is a subgroup of the automorphism group of $\G{P}$, as explained in~\S\ref{subsec:aut_P}. Since our methods make no assumptions on the codimension of $P$, by considering the automorphism group of $P\times\{1\}$ in $\Lambda\times\Z$ we are able to calculate the group of affine automorphisms $\AffAut{P}\leq\GL_n(\Z)\ltimes\Lambda$. As an illustration of our methods, we calculate the order of the automorphism group for each of the $473,\!800,\!776$ four-dimensional reflexive polytopes~\cite{KS00}: see Table~\ref{tab:num_4topes}.

The second approach -- to compute a normal form $\NF{P}$ for each isomorphism class -- is discussed in~\S\ref{sec:palp_normal_form}. This is the approach adopted by Kreuzer and Skarke in their {\palp} software~\cite{KS04}, and was used to construct the classification of three- and four-dimensional reflexive polytopes~\cite{KS98b,KS00}. Briefly, row and column permutations are applied to the vertex--facet pairing matrix $\PM$ of $P$, placing it in a form $\PMmax$ that is maximal with respect to a certain ordering. This in turn defines an order in which to list the vertices of $P$; the choice of basis is fixed by taking the Hermite normal form. In~\S\ref{subsec:affine_normal_form} we address how this can be modified to give an affine normal form for $P$, and in~\S\ref{subsec:palp_normal_form} we describe how {\palp} applies an additional reordering of the columns of $\PMmax$ before computing the normal form. The {\palp} source code for computing $\NF{P}$ is analyzed in detail in Appendix~\ref{apx:palp_source_code}.


In~\S\ref{sec:exploiting_aut} we address the problem of calculating $\PMmax$. We describe an inductive algorithm which attempts to exploit automorphisms of the matrix in order to simplify the calculation; pseudo-code is given in Appendix~\ref{apx:matrix_isomorphism}. Applying our algorithm to smooth Fano polytopes~\cite{Obr07}, which often have large numbers of symmetries, illustrates the advantage of this approach: see~\S\ref{subsec:analysis_smooth_db} and Table~\ref{tab:timings}. We end by giving, in~\S\ref{sec:laurent_normal_form}, an application of normal form to Laurent polynomials.

\subsection*{A note on implementation}
The algorithms described in~\S\ref{sec:iso_via_face_graph} were implemented using {\magma} in~$2008$ and officially released as part of {\magma}~V$2.16$~\cite{Magma,ConvChap}; {\palp} normal form was introduced by Kreuzer and Skarke in their {\palp} software~\cite{KS04} and reimplemented natively in {\magma}~V$2.18$ by the authors. The {\magma} algorithms\footnote{Users of {\magma} can freely view and edit the package code. The relevant files are contained in the subdirectory \texttt{package/Geometry/ToricGeom/polyhedron/}.}, including the reimplementation of {\palp} normal form, have recently been ported to the {\sage} project~\cite{sage} by Samuel Gonshaw\footnote{Gonshaw's implementation is available from \href{http://trac.sagemath.org/sage_trac/ticket/13525}{\texttt{http://trac.sagemath.org/sage\underline{\ }trac/ticket/13525}}.}, assisted by Tom Coates and the second author, and should appear in the~$5.6.0$ release.

\subsection*{Acknowledgments}
This work was motivated in part by discussions with Max Kreuzer during August and September 2010, shortly before his death that November. We are honoured that he found the time and energy for these conversations during this period. It forms part of the collaborative {\nolinebreak[4]{\palp}\small$++$\normalsize} project envisioned in~\cite{palp++}.

Our thanks to Tom Coates for many useful discussions, to Harald Skarke and Dmitrii Pasechnik for several helpful comments on a draft of this paper, to John Cannon for providing copies of the computational algebra software {\magma}, and to Andy Thomas for technical assistance. The first author was funded by a Summer Studentship as part of Tom Coates' Royal Society University Research Fellowship. The second author is supported by EPSRC grant EP/I008128/1.

\section{Isomorphism testing via the face graph}\label{sec:iso_via_face_graph}

\subsection*{Conventions}
Throughout this section we work with very general convex polytopes; we assume only that $P\subset\Lambda_\Q:=\Lambda\otimes\Q$ is a (non-empty) rational convex polytope, not necessarily of maximum dimension in the ambient lattice $\Lambda$.  The dual lattice $\Hom(\Lambda,\Z)$ is denoted by $\Lambda^*$.

Given two polytopes $P$ and $P'$, how can we decide whether they are isomorphic and, if they are, how can we construct an isomorphism between them? There are, of course, some obvious checks that can quickly provide a negative answer. We give a few examples, although this list is far from comprehensive.

\begin{itemize}
\item Do the dimensions of the polytopes agree?
\item Does $P$ contain the origin in its relative interior? Is the same true for $P'$?
\item Are $P$ and $P'$ both lattice polytopes?
\item Are the $f$-vectors of $P$ and $P'$ equal?
\item Do $P$ and $P'$ have the same number of primitive vertices?
\item Are $P$ and $P'$ simplicial? Are they simple?
\item If $P$ is of codimension one then there exists a unique hyperplane $H\subset\Lambda_\Q$ containing $P$, where $H=\{v\in\Lambda_\Q\mid\pro{v}{u}=k\}$ for some non-negative rational value $k$ and primitive dual lattice point $u\in\Lambda^*$. In particular, $k$ is invariant under change of basis. Does $k$ agree for both $P$ and $P'$?
\item If $P$ is of maximum dimension, any facet $F$ can be expressed in the form $F=\{v\in P\mid\pro{v}{u_F}=-c_F\}$, where $u_F\in\Lambda^*$ is a primitive inward-pointing vector normal to $F$, and $c_F\in\Q$ is the lattice height of $F$ over the origin. The value of $c_F$ is invariant under change of basis. Do the facet heights of $P$ and $P'$ agree, up to permutation?
\item If $P$ is a rational polytope, let $r_P$ be the smallest positive integer such that the dilation $r_PP$ is a lattice polytope. Do $r_P$ and $r_{P'}$ agree?
\end{itemize}

\begin{remark}
From a computational point of view, the intention with the above list is to suggest tests that are easy to perform. We assume that data such as the vertices and supporting hyperplanes of $P$ have already been calculated. Some computations, such as finding the $f$-vector, are more involved, but since the calculations will be required in what follows it seems sensible to use them at this stage.

In practice a number of other invariants may already be cached and could also be used: the volume $\Vol{P}$ or boundary volume $\Vol{\partial P}$; the number of lattice points $\abs{P\cap\Lambda}$ or boundary lattice points $\abs{\partial P\cap\Lambda}$; the Ehrhart $\delta$-vector; information about the polar polyhedron $P^*$. In particular cases some of this additional data may be easy to calculate; in general they are usually more time-consuming to compute than the isomorphism test described below.
\end{remark}

\begin{remark}
There are a few potential catches for the unwary when considering rational polytopes with $\dim{P}<\dim{\Lambda}$. For example, care needs to be taken when defining the supporting hyperplanes. Also, the notion of (normalised) volume $\Vol{P}$ requires some attention: the affine sublattice $\aff{P}\cap\Lambda$ may be empty, forcing us to either accept that $\Vol{P}$ can be undefined, or to employ interpolation. There is a natural dichotomy between those polytopes whose affine span contains the origin and those where $0\notin\aff{P}$. In the latter case, it is often better to consider the cone $C_P:=\cone{P}$ equipped with an appropriate grading such that dilations of $P$ can be realised by taking successive slices through $C_P$.
\end{remark}
\subsection{The labelled face graph}\label{subsec:labelled_face_graph}
In order to determine isomorphism we make use of the face graph $G(P)$ of $P$.
\begin{definition}
Let $P$ be an $n$-dimensional polytope with $f$-vector $(f_{-1},f_0,\ldots,f_n)$, where $f_k$ denotes the number of $k$-faces of $P$. By convention we set $f_{-1}=f_n=1$, representing, respectively, the empty set $\varnothing$ and the polytope $P$. The \emph{face graph} $G(P)$ is the graph consisting of $f_{-1}+f_0+\ldots+f_n$ vertices, where each vertex $v$ corresponds to a face $F_v$. Two vertices $v$ and $v'$ are connected by an edge if and only if $F_{v'}\subset F_v$ and $\dim{F_{v'}}=\dim{F_v}+1$.  Here the dimension of the empty face $\varnothing$ is taken to be ${-1}$.
\end{definition}

The face graph of a polytope is completely determined by the vertex--facet relations, and is the standard tool for determining combinatorial isomorphism of polytopes. We augment $G(P)$ by assigning labels to the vertices determined by some invariants of the corresponding face. Reducing a symmetry problem to the study of a (labelled) graph is a well-established computational technique: see, for example,~\cite{KaSc03,Pug05,MBW09,BSPRS12}. The intention is to decorate the graph with data capturing how $P$ lies in the underlying lattice $\Lambda$. To that end, we make the following definition.

\begin{definition}
For a point $v\in\Lambda_\Q$, let $u\in\Lambda$ be the unique primitive lattice point such that $v=\lambda u$ for some non-negative value $\lambda$ (set $u=0$, $\lambda=0$ if $v=0$). We define $\tilde{v}$ to be given by $\lceil\lambda\rceil u$, i.e. $\tilde{v}$ is the first lattice point after or equal to $v$ on the ray defined by $v$. Let $P\subset\Lambda_\Q$ be a polytope with vertices $\V{P}$. Then the \emph{index} $\idx{P}$ of $P$ is the index of the sublattice generated by $\left\{\tilde{u}\mid u\in\V{P}\right\}$ in $\linspan{P}\cap\Lambda$.
\end{definition}

\begin{definition}\label{def:labelled_face_graph}
Let $P$ be an $n$-dimensional polytope with face graph $G(P)$. To each vertex $v$ of $G(P)$ we assign the label
$$
\left\{\begin{array}{rl}
(\dim{F_v}),&\text{ if }F_v=\varnothing\text{ or }F_v=P;\\
(\dim{F_v},\idx{F_v}),&\text{ otherwise.}\\
\end{array}\right.
$$
We denote this labelled graph by $\G{P}$.
\end{definition}

\begin{remark}
In place of the index $\idx{F_v}$, it is tempting to use the volume $\Vol{F_v}$. However, computing the index is basic linear algebra, whereas computing the volume is generally difficult.
\end{remark}

\subsection{Additional labels}
When $P$ contains the origin strictly in its interior, we can make use of the \emph{special facets}. Recall from~\cite[\S3.1]{Obr07} that a facet $F$ is said to be special if $u\in\cone{F}$, where $u:=\sum_{v\in\V{P}}v$ is the sum of the vertices of $P$. Since $P$ contains the origin, there exists at least one special facet; we can extend the labelling to indicate which vertices of $\G{P}$ correspond to a special facet.

\begin{example}
The polytope $P:=\sconv{(1,0),(0,1),(-2,-3)}$ and its labelled face graph $\G{P}$ are depicted below. In the graph, the top-most vertex represents $P$ and the bottom vertex $\varnothing$. The sum of the vertices is $(-1,-2)$, so there is a unique special facet: the edge joining vertices $(1,0)$ and $(-2,-3)$ of index three, labelled $(1,3,1)$ in $\G{P}$. The edge joining vertices $(1,0)$ and $(0,1)$ is of index one and labelled $(1,1,0)$; the remaining edge is of index two and labelled $(1,2,0)$. The final entry of each facet label is used to indicate whether this is a special facet.

\centering
\raisebox{57pt}{\includegraphics[scale=1.1]{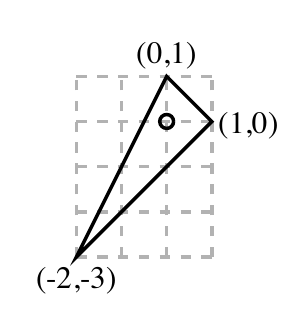}}\qquad
\includegraphics[scale=0.8]{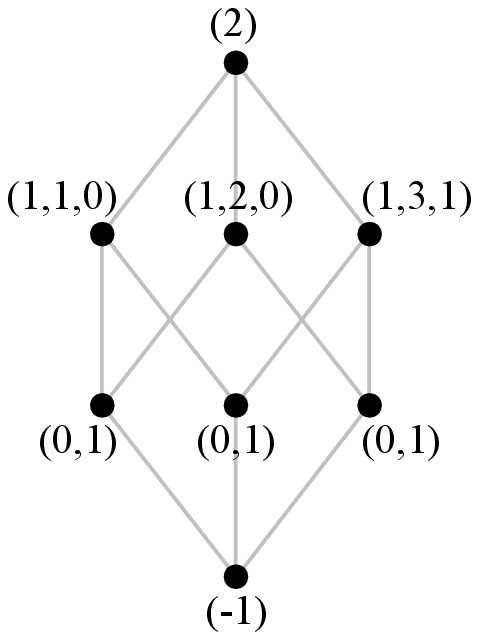}
\end{example}

If $P$ is a rational polytope, the vertices of $P$ provide an augmentation to the labelling of $\G{P}$. For any vertex $v\in\V{P}$ there exists a primitive lattice point $u\in\Lambda$ and a non-negative rational value $\lambda$ such that $v=\lambda u$. Since $\lambda$ is invariant under change of basis, the corresponding labels can be extended with this information. (Note that $\lambda=\idx{v}$ when $v$ is a lattice point, so this only provides additional information in the rational case.)

We do not claim that these are the only easily-computed invariants that can be associated with $\G{P}$. Other possibilities include encoding the linear relations between the vertices $\V{P}$ of $P$ in the graph labelling, and, in the maximum dimensional case, adding information about the lattice height of the supporting hyperplanes for each face.

\subsection{Recovering the isomorphism}
We now describe our algorithm for computing an isomorphism between two polytopes $P$ and $P'$. The initial step is to normalise the polytopes. If $P$ and $P'$ are not of maximum dimension in the ambient lattice $\Lambda$, then we first restrict to the sublattice $\linspan{P}\cap\Lambda$ (and, respectively, $\linspan{P'}\cap\Lambda$). It is possible that, even after restriction, $P$ and $P'$ are of codimension one. In that case, we work with the convex hull $\conv{P\cup\{0\}}$ (and similarly for $P'$). The important observations are that, after normalisation, $P$ is of maximum dimension, and that there exists at least one facet $F_0$ of $P$ such that $0\notin\aff{F_0}$.

Now we calculate an arbitrary graph isomorphism $\phi:\G{P}\rightarrow\G{P'}$. By restricting to the vertices of $\G{P}$ corresponding to the vertices $\V{P}$ of $P$, $\phi$ induces a map from the vertices of $P$ to the vertices of $P'$. The two polytopes $P$ and $P'$ are isomorphic only if $\phi$ exists, and any isomorphism $\Phi:\Lambda\rightarrow\Lambda$ mapping $P$ to $P'$ can be factored as $\phi\circ\chi$, where $\chi\in\Aut{\G{P}}$.

It remains to decide whether a particular choice of $\chi\in\Aut{\G{P}}$ determines a lattice isomorphism $\phi\circ\chi:\Lambda\rightarrow\Lambda$ sending $P$ to $P'$. For this we make use of the facet $F_0$. By construction $F_0$ is of codimension one, and does not lie in a hyperplane containing the origin. Hence there exists a choice of vertices $v_1,\ldots,v_n$ of $F_0$ which generate $\Lambda_\Q$ (over $\Q$). Denote the image of $v_i$ in $P'$ by $v'_i$, and consider the $n\times n$ matrices $V$ and $V'$ whose rows are given by, respectively, the $v_i$ and the $v'_i$. In order for $\phi\circ\chi$ to be a lattice map we require that $B:=V^{-1}V'\in\GL_n(\Z)$. In order for this to be an isomorphism from $P$ to $P'$ we require that $\{vB\mid v\in\V{P}\}=\V{P'}$.

\begin{remark}
We make two brief observations. First, in practice the automorphism group $\Aut{\G{P}}$ is often small. Second, it is an easy exercise in linear algebra to undo our normalisation process, lifting $B$ back to act on the original polytope.
\end{remark}

\subsection{Testing for equivalence}\label{subsec:equivalence}
Recall that two polytopes $P,P'\subset\Lambda_\Q$ are said to be \emph{equivalent} if there exists an isomorphism $B\in\GL_n(\Z)$ and a translation $c\in\Lambda$ such that $PB + c=P'$.

\begin{definition}
Let $\V{P}$ be the set of vertices of a polytope $P$. Then the \emph{vertex average} of $P$ is the point
$$b_P:=\frac{1}{\abs{\V{P}}}\sum_{v\in\V{P}}v\in\Lambda_\Q.$$
\end{definition}

Two polytopes $P$ and $P'$ are equivalent if and only if $b_P-b_{P'}\in\Lambda$ and $P-b_P$ is isomorphic to $P'-b_{P'}$.

\begin{example}
Consider the simplices
\begin{align*}
P&:=\sconv{(0,0,0),(2,1,1),(1,2,1),(1,1,2)},\\
P'&:=\sconv{(0,1,2),(1,0,0),(3,1,4),(4,2,6)}.
\end{align*}
The vertex averages are $b_P=(1,1,1)$ and $b_{P'}=(2,1,3)$, and $(P - b_P)B = P' - b_{P'}$, where
$$
B:=\small\begin{pmatrix}
2&1&3\\
-2&0&-1\\
1&0&1
\end{pmatrix}\normalsize
$$
Hence $P$ and $P'$ are equivalent.
\end{example}

\subsection{Determining the automorphism group of a polytope}\label{subsec:aut_P}
We can use the labelled face graph $\G{P}$ to compute the automorphism group $\Aut{P}$. We simply use the elements $\chi$ of $\Aut{\G{P}}$ to construct $\Aut{P}\le\GL_n(\Z)$. Notice that there is no requirement that $P$ is of maximum dimension in the ambient lattice $\Lambda$. Given this, we can also compute the affine automorphism group $\AffAut{P}$. Begin by embedding $P$ at height one in the lattice $\Lambda\times\Z$ (equivalently, consider the cone $C_P$ spanned by $P$ with appropriate grading). We refer to this embedded image of $P$ as $\tilde{P}$. The action of the automorphism group $\mathrm{Aut}\larger(\tilde{P}\larger)$ on $\tilde{P}$ restricts to an action on $P$, realising the full group of affine lattice automorphisms of~$P$. A detailed discussion of polyhedral symmetry groups and their applications can be found in~\cite{BEK84,BSS09,BSPRS12}.

\begin{example}\label{ex:involution}
Let $P$ be the three-dimensional simplicial polytope with seven vertices given by $(\pm1,0,0)$, $(0,\pm1,0)$, $(0,0,1)$, $(1,1,0)$, $(0,-1,-1)$. This is sketched below; the $f$-vector is $(1,7,15,10,1)$. The index $\idx{F}$ of each face $F$ is one (in fact $P$ is a smooth Fano polytope\footnote{Smooth Fano polytope number $13$ in the Graded Ring Database~\cite{GRDb}.}), and $P$ has four special facets (the four facets incident to the vertex $(1,0,0)$). 
The resulting labelled graph $\G{P}$ has automorphism group of order four, however $\Aut{P}$ has order two, and is generated by the involution $(0,0,1)\mapsto(0,-1,-1)$.

\centering
\vspace{10pt}
\includegraphics[scale=0.3]{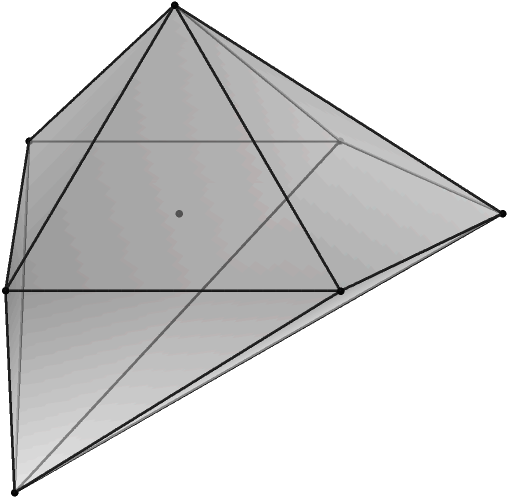}
\end{example}

\begin{example}\label{ex:24cell}
The four-dimensional centrally symmetric polytope $P$ with vertices
\begin{align*}
&\pm(1,0,0,0), \pm(0,1,0,0), \pm(0,0,1,0), \pm(0,0,0,1),\\
&\pm(1,-1,0,0), \pm(1,0,-1,0), \pm(1,0,0,-1), \pm(0,1,-1,0), \pm(0,1,0,-1),\\
&\pm(1,0,-1,-1), \pm(0,1,-1,-1), \pm(1,1,-1,-1)
\end{align*}
is the reflexive realisation of the $24$-cell, with $f$-vector $(1,24,96,96,24,1)$. It is unique amongst all $473,\!800,\!776$ reflexive polytopes in having $\abs{\Aut{P}}=1152$; in fact $\Aut{P}$ is isomorphic to the Weyl group $W(F_4)$. In particular, $P$ must be self-dual. The number of four-dimensional reflexive polytopes with $\abs{\Aut{P}}$ of given size are recorded in Table~\ref{tab:num_4topes}.
\end{example}

\begin{table}[htdp]
\centering
\begin{tabular}{cc}
\toprule
\begin{tabular}{cc}
$\abs{\Aut{P}}$&$\# P$\\
\midrule
\rowcolor[gray]{0.95}1&467705246\\
2&5925190\\
\rowcolor[gray]{0.95}3&1080\\
4&151416\\
\rowcolor[gray]{0.95}6&8218\\
8&6935\\
\rowcolor[gray]{0.95}10&4\\
12&1509\\
\rowcolor[gray]{0.95}16&756\\
18&2\\
\rowcolor[gray]{0.95}20&4\\
24&247\\
\rowcolor[gray]{0.95}32&23\\
\end{tabular}&
\begin{tabular}{cc}
$\abs{\Aut{P}}$&$\# P$\\
\midrule
\rowcolor[gray]{0.95}36&11\\
48&79\\
\rowcolor[gray]{0.95}64&5\\
72&10\\
\rowcolor[gray]{0.95}96&22\\
120&2\\
\rowcolor[gray]{0.95}128&2\\
144&2\\
\rowcolor[gray]{0.95}240&4\\
288&2\\
\rowcolor[gray]{0.95}384&6\\
1152&1\\
&\phantom{467705246}\\
\end{tabular}\\
\bottomrule
\end{tabular}
\vspace{0.5em}
\caption{The number $\# P$ of four-dimensional reflexive polytopes with automorphism group of size $\abs{\Aut{P}}$.}
\label{tab:num_4topes}
\end{table}

\begin{example}\label{ex:affaut}
Let $P=\sconv{(0,0),(1,0),(0,1)}$ be the empty simplex in $\Z^2$. Then $\Aut{P}$ is of order two, corresponding to reflection in the line $x=y$. To compute the affine automorphism group $\AffAut{P}$ of $P$, we consider $\tilde{P}=\sconv{(0,0,1),(1,0,1),(0,1,1)}$. The group $\mathrm{Aut}\larger(\tilde{P}\larger)$ is of order six, generated by
$$
\small\begin{pmatrix}-1&0&0\\-1&1&0\\1&0&1\end{pmatrix}
\normalsize\quad\text{and}\quad
\small\begin{pmatrix}0&-1&0\\1&-1&0\\0&1&1\end{pmatrix}.
$$
The first generator corresponds to the involution exchanging the vertices $(0,0)$ and $(1,0)$ of $P$, whilst the second generator corresponds to rotation of $P$ about its barycentre $(1/3,1/3)$, given by
$$
(x,y)\mapsto
(x-1/3,y-1/3)\begin{pmatrix}
0&-1\\
1&-1
\end{pmatrix}+(1/3,1/3)=(x,y)\begin{pmatrix}
0&-1\\
1&-1
\end{pmatrix}+(0,1).
$$
\end{example}

\section{Normal forms}\label{sec:palp_normal_form}
The method for determining isomorphism adopted by Kreuzer and Skarke in the software package {\palp}~\cite{KS04} is to generate a \emph{normal form} for the polytope $P$. We shall briefly sketch their approach. Their algorithm is described in detail in Appendix~\ref{apx:palp_source_code}.

Throughout we require that the polytope $P\subset\Lambda_\Q$ is a lattice polytope of maximum dimension. It is essential to the algorithm that the vertices are lattice points;  one could dilate a rational polytope by a sufficiently large factor to overcome this restriction, but in practice the resulting large vertex coefficients can cause computational problems of their own. Let $n$ denote the dimension of $P$, and $n_v$ be the number of vertices $\V{P}$. We can represent $P$ by an $n\times n_v$ matrix $V$ whose columns are given by the vertices. Obviously $V$ is uniquely defined only up to permutations $\sigma\in S_{n_v}$ of the columns.

Given any matrix $V$ with integer entries, we can compute its Hermite normal form $H(V)$. This has the property that, for all $B\in\GL_n(\Z)$, $H(V)=H(V\cdot B)$, however permuting the columns of $V$ will result in different Hermite normal forms. Na\"ively one could define the normal form $\NF{P}$ of $P$ to be
$$\min\left\{H(\sigma V)\mid\sigma\in S_{n_v}\right\},$$
where $\sigma V$ denotes the matrix obtained by permuting the columns of $V$ by $\sigma$, and the minimum is taken with respect to some ordering of the set of $n\times n_v$ integer matrices (say, lexicographic ordering). Unfortunately the size of $S_{n_v}$ is too large for this to be a practical algorithm.

\subsection{The pairing matrix}
The key to making this approach tractable is the vertex--facet pairing matrix.

\begin{definition}\label{defn:PM}
Let $P$ be a lattice polytope with vertices $v_j$, and let $(w_i,c_i)\in\Lambda^*\times\Z$ define the supporting hyperplanes of $P$; each $w_i$ is a primitive inward-pointing vector normal to the facet $F_i$ of $P$, such that $\pro{w_i}{v}=-c_i$ for all $v\in F_i$. The \emph{vertex--facet pairing matrix} $\PM$ is the $n_f\times n_v$ matrix with integer coefficients
$$\PM_{ij}:=\pro{w_i}{v_j}+c_i.$$
\end{definition}

In other words, the $ij$-th entry of $\PM$ correspond to the lattice height of $v_j$ above the facet $F_i$. This is clearly invariant under the action of $\GL_n(\Z)$. It is also invariant under (lattice) translation of $P$. Permuting the vertices of $P$ corresponds to permuting the columns of $\PM$, and permuting the facets of $P$ corresponds to permuting the rows of $\PM$. Thus there is an action of $S_{n_f}\times S_{n_v}$ on $\PM$: given $\sigma=(\sigma_f,\sigma_v)\in S_{n_f}\times S_{n_v}$,
$$(\sigma\PM)_{ij}:=\PM_{\sigma_f(i),\sigma_v(j)}.$$
There is a corresponding action on $V$ given by restriction:
$$\sigma V:=\sigma_v V.$$

Let $\PMmax$  denote the maximal matrix (ordered lexicographically) obtained from $\PM$ by the action of $S_{n_f}\times S_{n_v}$, realised by some element $\sigma_\text{max}$. Let $\Aut{\PMmax}\le S_{n_f}\times S_{n_v}$ be the automorphism group of $\PMmax$. Then:

\begin{definition}\label{defn:normal_form}
The \emph{normal form} of $P$ is
$$\NF{P}=\min\left\{H(\sigma\circ\sigma_\text{max} V)\mid\sigma\in\Aut{\PMmax}\right\}.$$
\end{definition}

\begin{remark}
Let $G$ be the group generated by the action of $\Aut{\PM}$ on the columns of $\PM$. Then $\Aut{P}\leq G$. Hence we have an alternative method for constructing the automorphism group when $P$ is a lattice polytope of maximum dimension.
\end{remark}

\begin{example}\label{eg:normal_form}
Consider the three-dimensional polytope $P$ with vertices $(1,0,0)$, $(0,1,0)$, $(0,0,1)$, $(-1,0,1)$, $(0,1,-1)$, $(0,-1,0)$, $(0,0,-1)$; $P$ is isomorphic to the polytope in Example~\ref{ex:involution} via the change of basis
$$\small\begin{pmatrix}
0&-1&-1\\
1&0&0\\
0&-1&0
\end{pmatrix}.$$
With the vertices in the order written above, and some choice of order for the facets, the vertex--facet pairing matrix is given by
$$\PM=\small\begin{pmatrix}
1&0&0&0&1&2&2\\
0&0&0&1&1&2&2\\
2&0&1&0&0&2&1\\
0&0&1&2&0&2&1\\
0&2&0&1&3&0&2\\
1&2&0&0&3&0&2\\
0&1&2&3&0&1&0\\
0&2&2&3&1&0&0\\
3&2&2&0&1&0&0\\
3&1&2&0&0&1&0
\end{pmatrix}\normalsize.\phantom{\PM=}$$
The maximum vertex--facet pairing matrix is
$$\PMmax=
\small\begin{pmatrix}
3&2&2&1&0&0&0\\
3&2&2&0&1&0&0\\
1&2&0&3&0&2&0\\
1&2&0&0&3&2&0\\
1&0&2&1&0&0&2\\
1&0&2&0&1&0&2\\
0&1&0&3&0&2&1\\
0&1&0&0&3&2&1\\
0&0&1&2&0&1&2\\
0&0&1&0&2&1&2\\
\end{pmatrix}\normalsize,\phantom{\PMmax=}$$
realised by, for example, the permutation $\left((1\ 5\ 2\ 6)(3\ 9)(4\ 10\ 7\ 8),(1\ 4\ 5)(3\ 6\ 7)\right)$ of $\PM$. The automorphism group of $\PMmax$ is of order two, generated by
$$\left((1\ 2)(3\ 4)(5\ 6)(7\ 8)(9\ 10),(4\ 5)\right).$$
We see that $\NF{P}$ is equal to
$$\small\begin{pmatrix}
1&0&1&0&-1&-1&0\\
0&1&-1&0&1&1&-1\\
0&0&0&1&-1&0&0\\
\end{pmatrix}\normalsize,$$
corresponding to the sequence of vertices $(1,0,0)$, $(0,1,0)$, $(1,-1,0)$, $(0,0,1)$, $(-1,1,-1)$, $(-1,1,0)$, and $(0,-1,0)$. In this example $\Aut{\PMmax}\cong\Aut{\NF{P}}$, and acts by exchanging the vertices $(0,0,1)$ and $(-1,1,-1)$.
\end{example}

\begin{example}\label{eg:matrix_aut_ne_poly_aut}
Let $P:=\sconv{(-1,-2,-2),(1,0,0),(0,2,1),(0,0,1)}$ be a three-dimensional reflexive polytope. This has
$$\phantom{,}\PMmax=\small\begin{pmatrix}
4&0&0&0\\
0&4&0&0\\
0&0&4&0\\
0&0&0&4\\
\end{pmatrix}\normalsize,\phantom{\PMmax=}$$
with $\Aut{\PMmax}\cong S_4$ of order $24$. However, $\abs{\Aut{P}}=8$; with ordering as above, the action on the vertices is given by the permutation group with generators  $(1\ 4\ 2\ 3)$ and $(3\ 4)$.
\end{example}

\subsection{Lattice polytopes of non-zero codimension}
Suppose that $P$ is a lattice polytope such that $\dim{P}<\dim{\Lambda}$. We can still define a normal form: how we proceed depends on whether $0\in\aff{P}$.

First suppose that $0\in\aff{P}$, so that $\aff{P}=\linspan{P}$. Set $d=\dim{P}$. We restrict $P$ to the sublattice $\linspan{P}\cap\Lambda\cong\Z^d$ and calculate the normal form there. The result can be embeded back into $\Lambda$ via
$$(a_1,\ldots,a_d)\mapsto(0,\ldots,0,a_1,\ldots,a_d).$$
Now suppose that $0\not\in\aff{P}$. In this case we consider the polytope $P_0:=\conv{P\cup\{0\}}$. The normal form $\NF{P_0}$ can be calculated and then the origin discarded.

\begin{example}
Let $P:=\sconv{(-1,1,1,0),(1,1,1,1),(0,0,0,-1)}$ be a lattice polygon of codimension two. The three-dimensional sublattice $\linspan{P_0}\cap\Z^4$ has generators $(1,0,0,0)$, $(0,1,1,0)$, and $(0,0,0,1)$. Let $\varphi:\Z^3\rightarrow\Z^4$ be the embedding given by right multiplication by the matrix
$$\small\begin{pmatrix}
1&0&0&0\\
0&1&1&0\\
0&0&0&1\\
\end{pmatrix}\normalsize.$$
Then $\varphi^*P_0$ has vertices $(-1,1,0)$, $(1,1,1)$, $(0,0,-1)$, and $(0,0,0)$, with normal form given by $(0,0,0)$, $(1,0,0)$, $(0,1,0)$, and $(1,1,2)$. Hence $\NF{P}$ corresponds to the vertices $(0,1,0,0)$, $(0,0,1,0)$, and $(0,1,1,2)$. In fact $P$ is isomorphic to $\NF{P}$ via the change of basis
$$\phantom{\in\GL_4(\Z).}\small\begin{pmatrix}
0&0&1&1\\
1&1&1&1\\
-1&0&0&0\\
0&-1&-1&-2\\
\end{pmatrix}\normalsize\in\GL_4(\Z).$$
\end{example}

\subsection{Affine normal form}\label{subsec:affine_normal_form}
The normal form can be adapted to give an \emph{affine normal form} $\AffNF{P}$ such that $\AffNF{P}=\AffNF{P'}$ if and only if polytopes $P$ and $P'$ are equivalent. One could simply define
$$\AffNF{P}:=\min\left\{\NF{P-v}\mid v\in\V{P}\right\}.$$
However, since the relative height of a vertex over a facet is unchanged by lattice translation, we have that $\PMmax$ is invariant. Hence
$$\AffNF{P}=\min\left\{H\mleft(\sigma\circ\sigma_\text{max} (V-v)\mright)\mid\sigma\in\Aut{\PMmax}, v\in\V{P}\right\}.$$

\begin{example}\label{eg:affine_normal_form}
Returning to the polytope in Example~\ref{eg:normal_form} we obtain
$$\phantom{.}\AffNF{P}=\small\begin{pmatrix}
0&1&0&0&3&2&1\\
0&0&1&0&2&1&2\\
0&0&0&1&-1&0&0\\
\end{pmatrix}\normalsize.\phantom{\AffNF{P}=}$$
\end{example}

\subsection{The {\palp} normal form}\label{subsec:palp_normal_form}
Kreuzer and Skarke's {\palp} normal form applies an additional modification to the order of the columns of the maximum vertex--facet pairing matrix $\PMmax$. For any $n_f\times n_v$ matrix $M$, let $c_M(j):=\max\left\{M_{ij}\mid 1\leq i\leq n_f\right\}$, and $s_M(j):=\sum_{i=1}^{n_f}M_{ij}$, where $1\leq j\leq n_v$. The following pseudo-code describes how the columns of $\PMmax$ (or, equivalently, the vertices of $P$) are rearranged.

\begin{algorithmic}
\State $M\gets\PMmax$
\For{$i=1$ \ForTo $n_v$}
\State $k\gets i$
\For{$j=i+1$ \ForTo $n_v$}
\If{$c_M(j)<c_M(k)\vee(c_M(j)=c_M(k)\wedge s_M(j)<s_M(k))$}
\State $k\gets j$
\EndIf
\EndFor
\State $M\gets\mathrm{SwapColumn}(M,i,k)$
\EndFor
\end{algorithmic}

\begin{example}\label{eg:palp_normal_form}
We revisit Example~\ref{eg:normal_form}. In this case, $\PMmax$ is modified by applying the permutation $(1\ 6\ 3\ 2)(4\ 7)$ to the columns, giving
$$\phantom{.}\small\begin{pmatrix}
2&2&0&0&0&3&1\\
2&2&0&0&1&3&0\\
2&0&2&0&0&1&3\\
2&0&2&0&3&1&0\\
0&2&0&2&0&1&1\\
0&2&0&2&1&1&0\\
1&0&2&1&0&0&3\\
1&0&2&1&3&0&0\\
0&1&1&2&0&0&2\\
0&1&1&2&2&0&0\\
\end{pmatrix}\normalsize.$$
The resulting {\palp} normal form corresponds to the sequence of vertices $(1,0,0)$, $(0,1,0)$, $(0,-1,0)$, $(-1,0,0)$, $(0,0,1)$, $(1,1,0)$, and $(0,-1,-1)$.
\end{example}

\begin{example}\label{eg:palp_affine_normal_form}
The affine normal form for the polytope in Example~\ref{eg:normal_form} with modified $\PMmax$ is given by
$$\phantom{.}\AffNF{P}=\small\begin{pmatrix}
0&1&1&2&0&0&2\\
0&0&2&2&0&-1&3\\
0&0&0&0&1&0&-1\\
\end{pmatrix}\normalsize.\phantom{\AffNF{P}=}$$
\end{example}

\section{Exploiting the automorphism group of the pairing matrix}\label{sec:exploiting_aut}
A crucial part of the normal form algorithm described in~\S\ref{sec:palp_normal_form} is the ability to efficiently calculate the maximum vertex--facet pairing matrix $\PMmax$. One also needs to know a permutation $\sigma$ such that $\sigma\PM=\PMmax$, and to be able to calculate $\Aut{\PMmax}$. These data can be constructed as $\PMmax$ is calculated -- this is the approach taken by the {\palp} source code described in Appendix~\ref{apx:palp_source_code} -- or recovered later. This section focuses on this second approach. A detailed algorithm is given in Appendix~\ref{apx:matrix_isomorphism}.

Consider a case when $\PM$ is very symmetric, so that the order of $\Aut{\PM}$ is large (for example, the vertex--facet pairing matrix for the $n$-dimensional polytope associated with projective space $\Proj^n$ has $\abs{\Aut{\PM}}=(n+1)!$). In such situations, the {\palp} algorithm is highly inefficient. Whilst computing $\PMmax$ the symmetries are not taken into account, so the algorithm needlessly explores equivalent permutations. Intuitively, one should be able to improve on the {\palp} algorithm by exploiting the automorphism group of $\PM$.

Given an $n_r\times n_c$ matrix $\PM$ and a group of possible column permutations $S$ (initially set to $S_{n_c}$), one can inductively convert this into $\PMmax$ as follows:
\begin{enumerate}
\item If $n_r=1$ then $\PMmax=\max\left\{\sigma\PM\mid\sigma\in S\right\}$.
\item If $\abs{S}=1$ then no permutations of the columns of $\PM$ are possible, and $\PMmax$ is given by sorting the rows of $\PM$ in decreasing order.
\item Suppose now that $n_r>1$ and $\abs{S}>1$.
\begin{enumerate}
\item Let $\mmax{R}:=\max\left\{\sigma\PM_i\mid\sigma\in S, 1\leq i\leq n_r\right\}$ be the largest row in $\PM$, up to the action of $S$.
\item Set $S':=\{\sigma\in S\mid \sigma\mmax{R}=\mmax{R}\}$.
\item For each row $1\leq i\leq n_r$ such that there exists a permutation $\sigma\in S$ with $\sigma\PM_i=\mmax{R}$, consider the matrix $M_{(i)}$ obtained from $\sigma\PM$ by deleting the $i$-th row. If $M_{(i)}\cong M_{(j)}$ for some $j<i$, then skip this case. Otherwise let $\mmax{M_{(i)}}$ be the $(n_r-1)\times n_c$ matrix obtained by inductively applying this process with $\PM\gets M_{(i)}$ and $S\gets S'$.
\item Set $\mmax{M}$ to be the maximum of all such $\mmax{M_{(i)}}$. Then
$$\PMmax=\left(\begin{array}{cc}
\mmax{R}\\
\hline
\mmax{M}
\end{array}\right).$$
\end{enumerate}
\end{enumerate}

\subsection{Test case: the database of smooth Fano polytopes}\label{subsec:analysis_smooth_db}
The algorithm described in Appendix~\ref{apx:matrix_isomorphism}, which we shall hereafter refer to as~\textsc{Symm}, was implemented by the authors and compared against the {\palp} algorithm. As Examples~\ref{eg:symmetry_vs_palp} and~\ref{eg:palp_vs_symmetry} illustrate, the difference in run-time between the two approaches can be considerable. 

\begin{example}\label{eg:symmetry_vs_palp}
Let $P$ be the six-dimensional polytope\footnote{Smooth Fano polytope number $1930$ in the Graded Ring Database~\cite{GRDb}.} with $14$ vertices
\begin{align*}
&\pm(1,0,0,0,0,0), \pm(0,1,0,0,0,0), \pm(0,0,1,0,0,0), \pm(0,0,0,1,0,0),\\
&\pm(0,0,0,0,1,0), \pm(0,0,0,0,0,1), \pm(1,1,1,1,1,1).
\end{align*}
The automorphism group $\Aut{\PM}$ is of order $10,\!080$. On our test machine the~{\palp} algorithm took $512.88$ seconds, whereas the~\textsc{Symm} algorithm took only $5.83$ seconds.
\end{example}

\begin{example}\label{eg:palp_vs_symmetry}
Let $P$ be the six-dimensional polytope\footnote{Smooth Fano polytope number $1854$ in the Graded Ring Database~\cite{GRDb}.} with $12$ vertices
\begin{align*}
&(1,0,0,0,0,0), (0,1,0,0,0,0), (0,0,1,0,0,0), (0,0,0,1,0,0), (0,0,0,0,1,0),\\
&(0,0,0,0,0,1), (-1,-1,-1,1,1,1), (0,0,1,-1,0,0), (0,0,-1,0,0,0),\\
&(0,1,1,-1,-1,-1), (0,-1,-1,0,0,0), (0,0,0,0,-1,-1).
\end{align*}
The automorphism group $\Aut{\PM}$ is of order $16$; the~{\palp} algorithm took $0.55$ seconds whilst the~\textsc{Symm} algorithm took $4.30$ seconds.
\end{example}

Table~\ref{tab:timings} contains timing data comparing the~{\palp} algorithm with the~\textsc{Symm} algorithm. This data was collected by sampling polytopes from {\O}bro's classification of smooth Fano polytopes~\cite{Obr07}. For each smooth polytope $P$ selected, the calculation was performed for both $P$ and $P^*$. In small dimensions the number of polytopes, and the time required for the computations, is small enough that the entire classification can be used. It is important to emphasise that the smooth Fano polytopes are atypical in that they can be expected to have a large number of symmetries, and so favour~\textsc{Symm}. Experimental evidence suggests that the ratio $r:=\abs{\Aut{\PM}}/{n_v}$ is a good proxy for deciding between the two choices. When $r<1$ the~{\palp} algorithm often performs better, whereas larger values indicate~\textsc{Symm} should be used.

\begin{table}[htdp]
\centering
\begin{tabular}{crrcrcrc}
\toprule
&&\multicolumn{2}{c}{{\palp}}&\multicolumn{2}{c}{\textsc{Symm}}&\multicolumn{2}{c}{\textsc{Best}}\\
Dim.&$\# P$&\multicolumn{1}{c}{Total}&\multicolumn{1}{c}{Average}&\multicolumn{1}{c}{Total}&\multicolumn{1}{c}{Average}&\multicolumn{1}{c}{Total}&\multicolumn{1}{c}{Average}\\
\cmidrule(lr){1-1} \cmidrule(lr){2-2} \cmidrule(lr){3-4} \cmidrule(lr){5-6} \cmidrule(lr){7-8}
\rowcolor[gray]{0.95}4&248&6.28&0.03&4.48&0.02&3.41&0.01\\
5&1732&98.30&0.06&59.53&0.03&46.17&0.03\\
\rowcolor[gray]{0.95}6&15244&6148.45&0.40&1510.32&0.10&1214.25&0.08\\
7&150892&152279.91&1.01&45230.55&0.30&34818.32&0.23\\
\rowcolor[gray]{0.95}8&281629&611795.13&2.17&152902.73&0.54&111426.70&0.40\\
\bottomrule
\end{tabular}
\vspace{0.5em}
\caption{Timing data, in seconds, for the~{\palp} algorithm and for the~\textsc{Symm} algorithm. The best possible time if one could infallibly choose the faster of the two algorithms is recorded by~\textsc{Best}.}
\label{tab:timings}
\end{table}

\section{Applications to Laurent polynomials}\label{sec:laurent_normal_form}
Let $f\in\C[x_1^{\pm1},\ldots,x_n^{\pm1}]$ be a Laurent polynomial in $n$ variables, and let $P:=\Newt{f}$ denote the Newton polytope of $f$. We require throughout that $\dim{P}=n$, i.e.~that $P$ is of maximum dimension in the ambient lattice. An element $B\in\GL_n(\Z)$ corresponds to the invertible monomial transformation
\begin{equation}\label{eq:cob}
\begin{array}{r@{\hspace{2pt}}c@{\hspace{2pt}}l}
\varphi_B:(\C^*)^n&\rightarrow&(\C^*)^n\\
x_j&\mapsto&x_1^{B_{1j}}\cdots x_n^{B_{nj}},
\end{array}
\end{equation}
and $g=\varphi_B^*f$ is also a Laurent polynomial. In particular, $\Newt{g}=P\cdot B$. 

As when working with lattice polytopes, it can be advantageous to be able to present $f$ in a normal form with respect to transformations of type~\eqref{eq:cob}.

\begin{definition}\label{defn:laurent_ordering}
Given two Laurent polynomials $f$ and $g$ such that $\Newt{f}=\Newt{g}$, we define an order $\preceq$ on $f$ and $g$ as follows. Let $v_1<v_2<\ldots<v_k$ be the lattice points in $\Newt{f}$, listed in lexicographic order. To each point $v_i$ there exists a (possibly zero) coefficient $c_i$ of $x^{v_i}$ in $f$, and coefficient $d_i$ in $g$. Define $\coeffs{f}:=(c_1,c_2,\ldots,c_k)$. We write $f\preceq g$ if and only if $\coeffs{f}\leq\coeffs{g}$.
\end{definition}

\begin{remark}
Any Laurent polynomial $f$ determines a pair $(\coeffs{f},\Newt{f})$. Conversely, given any pair $(c,P)$, where $c\in\C^k$ and $P\subset\Lambda_\Q$ is a maximum-dimensional lattice polytope such that $k=\abs{P\cap\Lambda}$, we can associate a Laurent polynomial. If we insist that the $c_i$ associated with the vertices $\V{P}$ are non-zero then we have a one-to-one correspondence.
\end{remark}

\begin{definition}\label{defn:laurent_normal_form}
Let $f$ be a Laurent polynomial, and set $P:=\Newt{f}$. Let $B\in\GL_n(\Z)$ be such that $P\cdot B=\NF{P}$. The \emph{normal form} for $f$ is
$$\NF{f}:=\mathrm{min}_\preceq\left\{\varphi_A\circ\varphi_B(f)\mid A\in\Aut{\NF{P}}\right\}.$$
\end{definition}

\begin{example}\label{eg:laurent_normal_form}
Consider the Laurent polynomial
$$f=2x^2y+\frac{1}{x}+\frac{3}{xy}.$$
Then $\NF{P}$ has vertices $(1,0)$, $(0,1)$, and $(-1,-1)$, with corresponding transformation matrix
$$B=\small\begin{pmatrix}
0&-1\\-1&1
\end{pmatrix}\normalsize\in\GL_2(\Z).$$
Under this transformation,
$$\varphi_B^*f=3x+y+\frac{2}{xy}$$
and $\coeffs{\varphi_B^*f}=(2,0,1,3)$. The automorphism group $\Aut{\NF{P}}\cong S_3$ acts by permuting the non-zero elements in the coefficient vector, hence
$$\NF{f}=3x+2y+\frac{1}{xy}.$$
\end{example}

A na\"ive implementation of Laurent normal form faces a serious problem: listing the points in a polytope is computationally expensive, and will often be the slowest part of the algorithm by many orders of magnitude. With a little care this can be avoided. What is really needed in Definition~\ref{defn:laurent_normal_form} is not the entire coefficient vector, but the closure of the non-zero coefficients under the action of $\Aut{\NF{P}}$. We illustrate this observation with an example.

\begin{example}\label{eg:orbit_closure}
Consider the Laurent polynomial
$$f=x^{50}y^{50}z^{50} + x^{50}y^{30} + \frac{x^{30}z^{30}}{y^{40}} + \frac{x^{10}}{y^{40}z^{20}} + xyz + \frac{y^{40}z^{20}}{x^{10}} + \frac{y^{40}}{x^{30}z^{30}} + \frac{1}{x^{50}y^{30}} + \frac{1}{x^{50}y^{50}z^{50}}.$$
Set $P=\Newt{f}$. Notice that $\abs{P\cap\Lambda}=285241$; enumerating the points in $P$ is clearly not the correct approach. The normal form $\NF{P}$ is given by change of basis
$$B=\small\begin{pmatrix}
-3&-4&-6\\
5&7&10\\
-12&-16&-23\\
\end{pmatrix}\normalsize\in\GL_3(\Z),$$
with
\begin{align*}
g:=\varphi_B^*f=x^{650}y^{880}z^{1270} + x^{500}y^{650}z^{950} + x^{10} + y^{10} + &\frac{1}{y^{10}} +\\
\frac{1}{x^{10}} + \frac{1}{x^{10}y^{13}z^{19}} + &\frac{1}{x^{500}y^{650}z^{950}} + \frac{1}{x^{650}y^{880}z^{1270}}.
\end{align*}
The automorphism group $G:=\Aut{\NF{P}}$ is of order two, generated by the involution $u\mapsto -u$. We consider the closure of the nine lattice points corresponding to the exponents of $g$ under the action of $G$. The only additional point is $(10,13,19)$. Thus we can express $\coeffs{g}$ with respect to these ten points:
$$\coeffs{g}=(1,1,1,1,1,1,1,0,1,1).$$
The key observation is that the action of $G$ on $g$ will not introduce any additional points, hence the lexicographically smallest coefficient sequence with respect to these points will also be the smallest coefficient sequence with respect to all the points of $\NF{P}$. By applying the involution we obtain the smaller coefficient sequence $(1,1,0,1,1,1,1,1,1,1)$, hence
\begin{align*}
\NF{f}=x^{650}y^{880}z^{1270} + x^{500}y^{650}z^{950} + x^{10}y^{13}z^{19} + &x^{10} + y^{10} +\\ 
\frac{1}{y^{10}} + \frac{1}{x^{10}} + &\frac{1}{x^{500}y^{650}z^{950}} + \frac{1}{x^{650}y^{880}z^{1270}}.
\end{align*}
\end{example}

We conclude this section by remarking that the automorphism group $\Aut{f}\leq\GL_n(\Z)$ of a Laurent polynomial $f$ can easily be constructed from $\Aut{\Newt{f}}$ by restricting to the subgroup that leaves $\coeffs{f}$ invariant.

\appendix
\section{The Kreuzer--Skarke algorithm}\label{apx:palp_source_code}
We describe in detail the algorithm used by Kreuzer and Skarke in {\palp}~\cite{KS04} to compute the normal form of a lattice polytope $P$ of maximum dimension $n$. Any such polytope can be represented by a $n\times n_v$ matrix $V$ whose columns correspond to the vertices of $P$. This matrix is unique up to permutation of columns and the action of $\GL_n(\Z)$; i.e. one can change the order of the vertices and the underlying basis for the lattice to obtain a different matrix $V'$.

The {\palp} normal form is a unique representation of the polytope $P$ such that if $Q$ is any other maximum dimensional lattice polytope, then $P$ and $Q$ are isomorphic if and only if their normal forms are equal. For any matrix $V$ with integer entries, and any $G\in\GL_n(\Z)$, the Hermite normal form of $G\cdot V$ is uniquely defined. The question is how to define a canonical order for the vertices, since permuting the vertices will lead to a different Hermite normal form.

In what follows, the line numbers refer to the {\palp} source file~\texttt{Polynf.c}\footnote{{\palp}~$1.1$, updated November~$2$,~$2006$. \href{http://hep.itp.tuwien.ac.at/~kreuzer/CY/palp/palp-1.1.tar.gz}{\texttt{http://hep.itp.tuwien.ac.at/$\sim$kreuzer/CY/palp/palp-1.1.tar.gz}}}. We have chosen our notation to correspond as closely as possible to the source code. The algorithm will be described in eight stages:
\begin{enumerate}
\item[\ref{asubsec:pairing_matrix}.] The pairing matrix;
\item[\ref{asubsec:max_pairing_matrix}.] The maximal pairing matrix;
\item[\ref{asubsec:first_row}.] Constructing the first row;
\item[\ref{asubsec:aut_step_1}.] Computing the restricted automorphism group, step I;
\item[\ref{asubsec:kth_row}.] Constructing the $k$-th row;
\item[\ref{asubsec:update_perms}.] Updating the set of permutations;
\item[\ref{asubsec:aut_step_2}.] Computing the restricted automorphism group, step II;
\item[\ref{asubsec:normal_form}.] Computing the normal form of the polytope.
\end{enumerate}

\subsection{The pairing matrix}\label{asubsec:pairing_matrix}
We start by constructing the pairing matrix $\PM$.

\vspace{1em}
\begin{center}
\begin{tabularx}{0.9\textwidth}{rX}
\toprule
Line:&197 (\texttt{Init\_rVM\_VPM})\\
Input:&A list of vertices and a list of equations for the supporting hyperplanes.\\
Output:&The pairing matrix $\PM$.\\
\bottomrule
\end{tabularx}
\end{center}
\vspace{1em}

Let $\left\{v_i\right\}_{i=1}^{n_v}$ be the vertices of $P$, in some order, and $\sum_{j=1}^nw_{ij}x_j+c_i=0$, $i=1,\ldots,n_f$, be the equations of the supporting hyperplanes of $P$. Here $n_v$ is equal to the number of vertices of $P$, and $n_f$ is equal to the number of facets of $P$. The $w_i$ are the inward-pointing primitive facet normals, and the $c_i$ are necessarily integers. The pairing matrix $\PM$ is the $n_f\times n_v$ matrix
$$\PM_{ij}=\sum_{k=1}^nw_{ik}v_{jk}+c_i=\pro{w_{i}}{v_{j}}+c_i$$
with integral coefficients.

The order of the columns of $\PM$ corresponds to an order of the vertices of $P$, and the order of the rows of $\PM$ corresponds to an order of the facets of $P$. Let $\rho=(r,c)\in S_{n_f}\times S_{n_v}$ act on $\PM$ via
$$(\rho\PM)_{ij}=\PM_{r(i)c(j)}.$$

\subsection{The maximal pairing matrix}\label{asubsec:max_pairing_matrix}
Let $\PMmax$ denote the maximal lexicographic matrix (when reading row by row) obtained from $\PM$ by reordering rows and columns, so that
$$\PMmax:=\max\left\{\rho\PM\mid\rho\in S_{n_f}\times S_{n_v}\right\}.$$

It can happen that $\Aut{\PM}\leq S_{n_f}\times S_{n_v}$ is non-trivial, say $\abs{\Aut{\PM}}=n_s$. Then we have $n_s$ permutations $\left\{\rho_i\right\}_{i=1}^{n_s}$ such that $\rho_i\PM=\PMmax$, and $n_s$ corresponding orders for the vertices of the polytope. Our main task is to compute $\PMmax$ and $\left\{\rho_i\right\}_{i=1}^{n_s}$ from $\PM$. This will be done by induction on the rows of $\PMmax$.

\subsection{Constructing the first row}\label{asubsec:first_row}
We begin by constructing the first row of $\PMmax$.

\vspace{1em}
\begin{center}
\begin{tabularx}{0.9\textwidth}{rX}
\toprule
Line:&348 (\texttt{Aux\_vNF\_Line})\\
Input:&The paring matrix $\PM$.\\
Output:&An array of permutations giving the first row of $\PMmax$.\\
\bottomrule
\end{tabularx}
\end{center}
\vspace{1em}

Set $n_s=1$ and maximise the first row of $\PM$, i.e. find a permutation $c_1\in S_{n_v}$ such that $\PM_{1c_1(i)}\leq\PM_{1c_1(j)}$, $j\leq i$:


\begin{algorithmic}
\State $n_s\gets 1$
\State $(r_1,c_1)\gets (1_{S_{n_f}},1_{S_{n_v}})$
\For{$j=1$ \ForTo $n_v$}
\State $m\gets\IndexOfMax{\PM_{1i}\mid i\geq j}$
\If{$m>1$}
\State $c_1\gets c_1(j\, m+j-1)$
\EndIf
\EndFor
\State $b\gets\PM_1$
\end{algorithmic}

Suppose we have computed the first $k-1$ lines, $n_s$ of which could be chosen to be the first row of $\PMmax$ (i.e.~up to reordering of the facets they are maximal among other lines and equal to the reference line, denoted $b$). Then we have integers $1\leq k_i\leq k-1<n_f$ with corresponding permutations $\rho_i=(r_i,c_i)\in S_{n_f}\times S_{n_v}$, $i=1,\ldots,n_s$, and a reference line defined by $b:=\PM_{k_1}$ such that:
$$\PM_{k_ic_i(j)}=b_{c_1(j)},\qquad i=1,\ldots,n_s, j=1,\ldots,n_v.$$
Set $r_i=(1\, k_i)$ to be the permutation which moves the line in question to the first row of $\PM$. Now we consider the $k$-th row of $\PM$. Find the maximal element $\max_j\{\PM_{kj}\}$, say $\PM_{km}$, and let $c_{n_s+1}=(1\, m)$. We compare this against the reference line. If $\PM_{kc_{n_s+1}(1)}<b_{c_1(1)}$ then continue with the next line (or stop if we are at the last line), otherwise continue constructing $c_{n_s+1}$. If $\max_{j>1}\{\PM_{kc_{n_s+1}(j)}\}=\PM_{kc_{n_s+1}(m)}$ then let $c_{n_s+1}\mapsto c_{n_s+1}\,(2\, m)$ and verify that $\PM_{kc_{n_s+1}(2)}<b_{c_1(2)}$; if this inequality fails to hold then continue with the next element.

If the line $k$ is not less than the reference line $b$ then we set $r_{n_s+1}=(1\, k)$ and have two cases to consider:
\begin{enumerate}
\item If $\PM_{kc_{n_s+1}(j)}=b_{c_1(j)}$, $j=1,\ldots,n_v$, then we have a new case of symmetry. We set $k_{n_s+1}:=k$ and increment the number of symmetries $n_s$.
\item Otherwise we have found a (lexicographically) bigger row and so obtain a new reference line. We set $b:=\PM_k$, $k_1:=k$, and $\rho_1:=(r_{n_s+1},c_{n_s+1})$, and reset the number of symmetries $n_s$.
\end{enumerate}


\begin{algorithmic}
\For{$k=2$ \ForTo $n_f$}
\State $(r_{n_s+1},c_{n_s+1})\gets (1_{S_{n_f}},1_{S_{n_v}})$
\State $m\gets\IndexOfMax{\PM_{kc_{n_s+1}(j)}\mid j\geq 1}$
\If{$m>1$}
\State $c_{n_s+1}\gets c_{n_s+1}(1\, m)$
\EndIf
\State $d\gets\PM_{kc_{n_s+1}(1)}-b_{c_1(1)}$
\If{$d<0$}
\Continue
\EndIf
\For{$i=2$ \ForTo $n_v$}
\State $m\gets\IndexOfMax{\PM_{kc_{n_s+1}(j)}\mid j\geq i}$
\If{$m>1$}
\State $c_{n_s+1}\gets c_{n_s+1}(i\, m+i-1)$
\EndIf
\If{d=0}
\State $d\gets\PM_{kc_{n_s+1}(i)}-b_{c_1(i)}$
\If{$d<0$}
\Break
\EndIf
\EndIf
\EndFor
\If{$d<0$}
\Continue
\EndIf
\State $r_{n_s+1}\gets r_{n_s+1}(1\, k)$
\If{$d=0$}
\State $n_s\gets n_s+1$
\Else
\State $(r_1,c_1)\gets (r_{n_s+1},c_{n_s+1})$
\State $n_s\gets 1$
\State $b\gets\PM_k$
\EndIf
\EndFor
\end{algorithmic}

\subsection{Computing the restricted automorphism group, step I}\label{asubsec:aut_step_1}
Once the first row of $\PMmax$ has been constructed, it imposes restrictions on any future column permutations: they must fix the first row.

\vspace{1em}
\begin{center}
\begin{tabularx}{0.9\textwidth}{rX}
\toprule
Line:&376 (\texttt{Aux\_vNF\_Line})\\
Input:&The first line of the maximal pairing matrix.\\
Output:&The array $S$ capturing its automorphism group.\\
\bottomrule
\end{tabularx}
\end{center}
\vspace{1em}

Suppose that the row is equal to blocks of $a_i$'s, each of size $n_i$ , $i=1,\ldots,k$, where $\sum_{i=1}^k n_i=n_v$:
$$
\left(\begin{array}{ccc|ccc|c|ccc}
a_1&\ldots&a_1&a_2&\ldots&a_2&\ldots&a_k&\ldots&a_k
\end{array}\right).$$

It is clear that if we had such a row, the only permutations of columns allowed in the construction of later rows will be those factoring through $S_{n_1}\times S_{n_2}\times\ldots\times S_{n_k}$. The symmetry of this row is encoded in an array $S$ such that if $S(i)=j$ and $S(S(i))=S(j)=h$ then the index $i$ is in the block delimited by the indices $j$ and $h$ (depending on whichever is greater). We represent $S$ as an array
$$
\begin{array}{r}
\left(\begin{array}{cccc|cccc|c}
n_1&1&\ldots&1&n_1+n_2&n_1+1&\ldots&n_1+1&\ldots\phantom{xxxxxxxx}\end{array}\right.\\
\left.\begin{array}{|cccc}
n_v&1+\sum_{i=1}^{k-1}n_i&\ldots&1+\sum_{i=1}^{k-1}n_i\end{array}\right)
\end{array}
$$

\begin{example}\label{aex:first_S}
The symmetries of the row $\left(5\ 5\ 5\ 5\ 4\ 3\ 3\ 2\ 2\ 2\ 1\ 0\ 0\right)$ are encoded by the array
$$
S=\left(\begin{array}{cccc|c|cc|ccc|c|cc}
4&1&1&1&5&7&6&10&8&8&11&13&12
\end{array}\right).
$$
\end{example}

When $S=\left(1\ 2\ \ldots\ n_v\right)$ the columns are fixed and we may only permute the rows. The computation of $S$ is summarised in the following pseudo-code:


\begin{algorithmic}
\State $S\gets\left(1\ 2\ \ldots\ n_v\right)$
\For{$i=2$ \ForTo $n_v$}
\If{$\PM_{r_1(1)c_1(i-1)}=\PM_{r_1(1)c_1(i)}$}
\State $S(i)\gets S(i-1)$
\State $S^2(i)\gets S(S(i))+1$
\Else
\State $S(i)\gets i$
\EndIf
\EndFor
\end{algorithmic}

\subsection{Constructing the $k$-th row}\label{asubsec:kth_row}
Proceeding by induction on the rows, we construct the remaining rows of $\PMmax$.

\vspace{1em}
\begin{center}
\begin{tabularx}{0.9\textwidth}{rX}
\toprule
Line:&289 (\texttt{Aux\_vNF\_Line})\\
Input:&$\PM$, the permutations $\{p_i\}_{i=1}^{n_s}$, and the array $S$.\\
Output:&The $k$-th line of the maximal pairing matrix.\\
\bottomrule
\end{tabularx}
\end{center}
\vspace{1em}


Assume we have computed the first $l-1<n_f-1$ rows of $\PMmax$ and the associated symmetry array $S$ (notice that the last row of $\PMmax$ need not be computed as it is completely determined), together with $n_s$ distinct permutations $\rho_i=(r_i,c_i)\in S_{n_f}\times S_{n_v}$ such that
$$\mmax{\PM_{kj}}=\PM_{r_i(k)c_i(j)}\qquad\text{ for all }1\leq j\leq n_v, 1\leq k<l, 1\leq i\leq n_s.$$

We have to consider each configuration given by the permutations $\left\{\rho_i\right\}_{i=1}^{n_s}$. For each configuration we generally obtain $n_\rho$ ways to construct the line $l$, moreover some constructions might give a smaller line, hence $n_s$ will have to be updated as we proceed. Let $\tilde{n}_{s}$ record the initial value of $n_s$.

First consider the case $k=\tilde{n}_{s}$. We will construct a candidate line for the $l$-th row of $\PMmax$; this will be our reference line against which the other cases will be compared. If a greater candidate is found, all the preceding computations will have to be deleted and redone with the new candidate. If a given case lead to a smaller line than the reference, it will have to be deleted.

Initially set the local number of symmetries, $n_\rho$, to zero and initialise the permutation $\tilde{\rho}_{n_\rho}=\rho_k$. We start with the line $\tilde{r}_{n_\rho}(l)$ by finding the maximal element of the first symmetry block. Suppose that
$$\max\left\{\PM_{\tilde{r}_{n_\rho}(l)\tilde{c}_{n_\rho}(i)}\mid 1\leq i\leq S(1)\right\} = \PM_{\tilde{r}_{n_\rho}(l)\tilde{c}_{n_\rho}(m)}.$$
Then we update $\tilde{c}_{n_\rho}$ to $\tilde{c}_{n_\rho}\,(1\, m)$. This maximal value is saved in the reference line which we denote $l_r$ (if it were already constructed, $k<\tilde{n}_{s}$, we move straight to the tests below). We increment $n_\rho$ by one to reflect this new candidate, initialise $\tilde{\rho}_{n_\rho}=p_{k}$, and proceed to consider the maximal entries in the first symmetry block for other lines $r_{k}(s)$, $s=l+1,\ldots,n_f$.

Inductively, suppose we have considered $s-1$ lines where $n_\rho$ of them have a maximal element in the first symmetry block equal to the one of the reference line $l_r(1)$, and the others have smaller values. We also have $\tilde{r}_{n_\rho}=r_{k}$ from the initialisation. Consider the line $\tilde{r}_{n_\rho}(s)$ and find its maximal element in $1,\ldots,S(1)$ as above, updating $\tilde{c}_{n_\rho}$. Now if $\PM_{\tilde{r}_{n_\rho}(s)\tilde{c}_{n_\rho}(1)}<l_r(1)$ then proceed to the case $s+1$, if possible. Otherwise $\tilde{r}_{n_\rho}\mapsto\tilde{r}_{n_\rho}\,(l\, s)$ and there are two possibilities: if $\PM_{\tilde{r}_{n_\rho}(s)\tilde{c}_{n_\rho}(1)}=l_r(1)$ then increase the number of symmetries $n_\rho\mapsto n_\rho+1$ and move to $s+1$, after initialising the new permutation $\tilde{\rho}_{n_\rho}=\rho_{k}$; if $\PM_{\tilde{r}_{n_\rho+1}(s)\tilde{c}_{n_\rho+1}(1)}>l_r(1)$ then redefine the first element of the reference line $l_r(1):=\PM_{\tilde{r}_{n_\rho}(s)\tilde{c}_{n_\rho}(1)}$, update the first permutation $\tilde{\rho}_0=\tilde{\rho}_{n_\rho}$, and reset $n_\rho=1$ ready for the next permutation $\tilde{\rho}_{n_\rho}=\rho_k$.



\begin{algorithmic}
\State $c\gets 1$
\State $n_\rho\gets 0$
\State $ccf\gets cf$
\State $(\tilde{r}_{n_\rho},\tilde{c}_{n_\rho})\gets (r_k,c_k)$
\For{$s=l$ \ForTo $n_f$}
\For{$j=2$ \ForTo $S(1)$}
\If{$\PM_{\tilde{r}_{n_\rho}(s)\tilde{c}_{n_\rho}(c)}<\PM_{\tilde{r}_{n_\rho}(s)\tilde{c}_{n_\rho}(j)}$}
\State $\tilde{c}_{n_\rho}\gets\tilde{c}_{n_\rho}(c\, j)$
\EndIf
\EndFor
\If{$ccf=0$}
\State $l_r(1)\gets\PM_{\tilde{r}_{n_\rho}(s)\tilde{c}_{n_\rho}(1)}$
\State $\tilde{r}_{n_\rho}\gets\tilde{r}_{n_\rho}(l\, s)$
\State $n_\rho\gets n_\rho+1$
\State $ccf\gets 1$
\State $(\tilde{r}_{n_\rho},\tilde{c}_{n_\rho})\gets (r_k,c_k)$
\Else
\State $d\gets\PM_{\tilde{r}_{n_\rho}(s)\tilde{c}_{n_\rho}(1)}-l_r(1)$
\If{$d<0$}
\Continue
\ElsIf{$d=0$}
\State $\tilde{r}_{n_\rho}\gets\tilde{r}_{n_\rho}(l\, s)$
\State $n_\rho\gets n_\rho+1$
\State $(\tilde{r}_{n_\rho},\tilde{c}_{n_\rho})\gets (r_k,c_k)$
\Else
\State $l_r(1)\gets\PM_{\tilde{r}_{n_\rho}(s)\tilde{c}_{n_\rho}(1)}$
\State $cf\gets 0$
\State $(\tilde{r}_1,\tilde{c}_1)\gets (\tilde{r}_{n_\rho},\tilde{c}_{n_\rho})$
\State $n_\rho\gets 1$
\State $(\tilde{r}_{n_\rho},\tilde{c}_{n_\rho})\gets (r_k,c_k)$
\State $n_s\gets k$
\State $\tilde{r}_{n_\rho}\gets \tilde{r}_{n_\rho}(l\, s)$
\EndIf
\EndIf
\EndFor
\end{algorithmic}

Note that the initial value of the \emph{comparison flag} $cf$ is $0$. This indicates that the reference line has not been initialised; it is also reset to zero when a greater candidate is found. We will see later how $cf$ is updated.

We need to construct other elements of $l_r$. 
Inductively, suppose we are constructing the entry $i$ of $l_r$ and we have $n_\rho$ symmetries with corresponding permutations $\tilde{\rho}_{j}$, $j=0,\ldots,n_\rho-1$. If $n_\rho=0$ we move to the next configuration $k-1$ after having updated the symmetries accordingly, i.e.~we do not save the current configuration. Otherwise, start with the last $j=n_\rho-1$. Determine where the corresponding block of symmetry ends for $i$ by looking at the maximum of $S(i)$ and $S^{2}(i)$, which we will call $h$. Then compute
$$\max\left\{\PM_{\tilde{r}_{j}(l)\tilde{c}_{j}(\lambda)}\mid i\leq\lambda\leq h\right\} =\PM_{\tilde{r}_{j}(l)\tilde{c}_{j}(m)}$$
and update $\tilde{c}_{j}\mapsto\tilde{c}_{j}\,(i\, m)$ . This value is saved in the reference line $l_r(i)$.
Then we consider (inductively) any cases of symmetry with $j<n_\rho-1$ and compute the $i$-th entry in the same manner as above: if $\PM_{\tilde{r}_{j}(l)\tilde{c}_{j}(i)}=l_r(i)$ then continue with the next $j$; if $\PM_{\tilde{r}_{j}(l)\tilde{c}_{j}(i)}<l_r(i)$ then the current case is removed and we update $n_\rho\mapsto n_\rho-1$; finally if $\PM_{\tilde{r}_{j}(l)\tilde{c}_{j}(i)}>l_r(i)$ then all
cases previously considered are irrelevant, so we let $n_\rho=j+1$ and the reference line is updated $l_r(i)=\PM_{\tilde{r}_{j}(l)\tilde{c}_{j}(i)}$. 



\begin{algorithmic}
\For{$c=2$ \ForTo $n_v$}
\State $h\gets S(c)$
\State $ccf\gets cf$
\If{$h<c$}
\State $h\gets S(h)$
\EndIf
\State $s\gets n_\rho$
\While{$s>0$}
\State $s\gets s-1$
\For{$j=c+1$ \ForTo $h$}
\If{$\PM_{\tilde{r}_s(l)\tilde{c}_s(c)}<\PM_{\tilde{r}_s(l)\tilde{c}_s(j)}$}
\State $\tilde{c}_s\gets \tilde{c}_s(cj)$
\EndIf
\EndFor
\If{$ccf=0$}
\State $l_r(c)\gets\PM_{\tilde{r}_s(l)\tilde{c}_s(c)}$
\State $ccf\gets 1$
\Else
\State $d\gets\PM_{\tilde{r}_s(l)\tilde{c}_s(c)}-l_r(c)$
\If{$d<0$}
\State $n_\rho\gets n_\rho-1$
\If{$n_\rho>s$}
\State $(\tilde{r}_s,\tilde{c}_s)\gets (\tilde{r}_{n_\rho},\tilde{c}_{n_\rho})$
\EndIf
\ElsIf{$d>0$}
\State $l_r(c)\gets\PM_{\tilde{r}_s(l)\tilde{c}_s(c)}$
\State $cf\gets 0$
\State $n_\rho\gets s+1$
\State $n_s\gets k$
\EndIf
\EndIf
\EndWhile
\EndFor
\end{algorithmic}

\subsection{Updating the set of permutations}\label{asubsec:update_perms}
The last step in the construction of the line $l$ is to organise the new symmetries for a given case $k$.

\vspace{1em}
\begin{center}
\begin{tabularx}{0.9\textwidth}{rX}
\toprule
Line:&333 (\texttt{Aux\_vNF\_Line})\\
Input:&The permutations $\{p_i\}_{i=1}^{n_s}$ and the newly computed $\{\tilde{\rho}_i\}_{i=0}^{n_\rho-1}$.\\
Output:&The updated set $\{p_i\}_{i=1}^{n_s}$.\\
\bottomrule
\end{tabularx}
\end{center}
\vspace{1em}

Recall that $\tilde{n}_{s}$ denotes the number of symmetries we had before performing the computations for the line $l$ of $\PMmax$, and $n_s\leq\tilde{n}_s$ represents the updated number symmetries.
Our current construction of the line $l$ may well introduce new symmetries, 
so-called~\emph{local symmetries}, of which there are $n_\rho$. We can have $n_\rho=0$, in which case all the configurations in the case $k$ lead to a smaller candidate for $l$. When $n_\rho>0$ the local symmetries are represented by the set $\{\tilde{\rho}_i\}_{i=0}^{n_\rho-1}$ of new permutations.

We now update the array of all permutations. If $n_s>k$ we set $\rho_k=\rho_{n_s}$; we want the set of permutations $\{\rho_i\}_{i=1}^{n_s}$ to be updated so that the only cases which need to be considered are those with index $i<k$. Since we are appending $n_\rho$ new permutations at end for the indices $i\geq n_s$, so $n_s\rightarrow n_s+n_\rho-1$. If $n_\rho=0$ then nothing is appended and $n_s$ decreases by one as required. Finally, we update the comparison flag $cf$ to reflect the current number of symmetries.


\begin{algorithmic}
\State $n_s\gets n_s-1$
\If{$n_s>k-1$}
\State $(r_k,c_k)\gets (r_{n_s+1},c_{n_s+1})$
\EndIf
\State $cf\gets n_s+n_\rho$
\For{$s=0$ \ForTo $n_\rho-1$}
\State $(r_{n_s+1},c_{n_s+1})\gets (\tilde{r}_s,\tilde{c}_s)$
\State $n_s\gets n_s+1$
\EndFor
\end{algorithmic}

\subsection{Computing the restricted automorphism group, step II}\label{asubsec:aut_step_2}
Once a new row of $\PMmax$ has been compute we need to update $S$ to reflect the symmetries of this row. This is done by restricting the blocks previously delimited by $S$ to reflect any additional constraints imposed by the row.

\begin{example}
Continuing Example~\ref{aex:first_S}, suppose that the second row of the candidate $\PMmax$ has been computed, and that the two rows are given by
$$
\left(\begin{array}{ccccccccccccc}
5&5&5&5&4&3&3&2&2&2&1&0&0\\
4&3&3&3&3&2&2&2&1&0&0&0&0
\end{array}\right).
$$
The corresponding array $S$ is
$$
\left(\begin{array}{c|ccc|c|cc|c|c|c|c|cc}
1&4&2&2&5&7&6&8&9&10&11&13&12
\end{array}\right).
$$
\end{example}

\vspace{1em}
\begin{center}
\begin{tabularx}{0.9\textwidth}{rX}
\toprule
Line:&376 (\texttt{Aux\_vNF\_Line})\\
Input:&The newly computed upper block of the maximal pairing matrix.\\
Output:&The updated array $S$ capturing the automorphism group of the matrix.\\
\bottomrule
\end{tabularx}
\end{center}
\vspace{1em}



\begin{algorithmic}
\State $c\gets 1$
\While{$c<n_v+1$}
\State $s\gets S(c)+1$
\State $S(c)\gets c$
\State $c\gets c+1$
\While{$c<s$}
\If{$\PM_{r_1(l)c_1(c)}=\PM_{r_1(l)c_1(c-1)}$}
\State $S(c)\gets S(c-1)$
\State $S^2(c)\gets S^2(c)+1$
\Else
\State $S(c)\gets c$
\EndIf
\State $c\gets c+1$
\EndWhile
\EndWhile
\end{algorithmic}

\subsection{Computing the normal form of the polytope}\label{asubsec:normal_form}
Inductively, we have obtained $n_s$ permutations $\left\{\rho_i=(r_i,c_i)\right\}_{i=1}^{n_s}$ such that $\rho_i\PM=\PMmax$. We are really only interested in the permutations of the columns, since they correspond to permutations of the vertices of $P$. The {\palp} algorithm computes a new order for the columns of $\PMmax$ based on the following: the maximum coefficient in the column; the sum of the coefficients in the column; and the relative position of the column in $\PMmax$. Let $\rho_c\in S_{n_v}$ denote this column permutation.

\vspace{1em}
\begin{center}
\begin{tabularx}{0.9\textwidth}{rX}
\toprule
Line:&216 (\texttt{New\_pNF\_Order})\\
Input:&The maximal pairing matrix $\PMmax$.\\
Output:&The column permutation $\rho_c\in S_{n_v}$.\\
\bottomrule
\end{tabularx}
\end{center}
\vspace{1em}



\begin{algorithmic}
\State $\PMmax\gets p_1\PM$
\State $p_c\gets 1_{S_{n_v}}$
\State $\mmax{M}\gets\left\{\max_{1\leq i\leq n_f}\left\{\mmax{\PM_{ij}}\right\}\mid 1\leq j\leq n_v\right\}$ 
\State $\mmax{S}\gets\left\{\sum_{1\leq i\leq n_f}\mmax{\PM_{ij}}\mid 1\leq j\leq n_v\right\}$
\For{$i=1$ \ForTo $n_v$}
\State $k\gets i$
\For{$j=i+1$ \ForTo $n_v$}
\If{$(\mmax{M_j}<\mmax{M_k})\vee((\mmax{M_j}=\mmax{M_k})\wedge(\mmax{S_j}<\mmax{S_k}))$}
\State $k\gets j$
\EndIf
\EndFor
\If{$k\neq i$}
\State $\mmax{M}\gets\mathrm{SwapRow}(\mmax{M},i,k)$
\State $\mmax{S}\gets\mathrm{SwapRow}(\mmax{S},i,k)$
\State $p_c\gets p_c(i\, k)$
\EndIf
\EndFor
\end{algorithmic}


Given the column permutations $\rho_c$ and $c_i$, $i=1\ldots,n_s$, we obtain a permutation of the vertices of $P$, and hence of the columns of the vertex matrix $V$. We let $V_i$ denote this reordered vertex matrix. The remaining freedom -- the action of $\GL_n(\Z)$ corresponding to the choice of lattice basis -- is removed by computing the Hermite normal form $H(V_i)$.

\vspace{1em}
\begin{center}
\begin{tabularx}{0.9\textwidth}{rX}
\toprule
Line:&134 (\texttt{GLZ\_Make\_Trian\_NF})\\
Input:&A matrix with integer coefficents.\\
Output:&The Hermite normal form of the matrix.\\
\bottomrule
\end{tabularx}
\end{center}
\vspace{1em}


The {\palp} normal form is simply the minimum amongst the $H(V_i)$.

\vspace{1em}
\begin{center}
\begin{tabularx}{0.9\textwidth}{rX}
\toprule
Line:&399 (\texttt{Aux\_Make\_Triang})\\
Input:&The column permutations $\rho_c$ and $\{c_i\}_{i=1}^{n_s}$, and the vertex matrix $V$.\\
Output:&The normal form.\\
\bottomrule
\end{tabularx}
\end{center}
\vspace{1em}


\section{Calculating the maximum pairing matrix}\label{apx:matrix_isomorphism}
Let $M$ be an $n_r\times n_c$ matrix. Recall that we define an action of $\sigma=(\sigma_r,\sigma_c)\in S_{n_r}\times S_{n_c}$ on the rows and columns of $M$ via $(\sigma M)_{ij}:=M_{\sigma_r(i),\sigma_c(j)}$, and that we call two matrices $M$ and $M'$ isomorphic if there exists some permutation $\sigma\in S_{n_r}\times S_{n_c}$ such that $\sigma(M)=M'$. We begin by briefly describing one approach to determining when two matrices are isomorphic.

Given a matrix $M$, we associate a bipartite graph $G(M)$ with $n_r+n_c$ vertices, where the vertices $v_i$, $v_{n_r+j}$ are connected by an edge $E_{ij}$ for all $1\leq i\leq n_r$, $1\leq j\leq n_c$. Each edge $E_{ij}$ is labelled with the corresponding value $M_{ij}$. The vertices $v_i$, $1\leq i\leq n_r$, are labelled with one colour, whilst the vertices $v_{n_r+j}$, $1\leq j\leq n_c$, are labelled with a second colour. This distinguishes between vertices representing rows of $M$ and vertices representing columns of $M$. Clearly two matrices $M$ and $M'$ are isomorphic if and only if the graphs $G(M)$ and $G(M')$ are isomorphic. We note also that the automorphism group $\Aut{M}\leq S_{n_r}\times S_{n_c}$ is given by the automorphism group of $G(M)$.

We now describe a recursive algorithm to compute $\PMmax$ from $\PM$. For readability, we shall split this algorithm into three parts, with a brief discussion preceding each part.

\vspace{1em}
\begin{center}
\begin{tabularx}{0.9\textwidth}{rX}
\toprule
Input:&A matrix $\PM$.\\
Output:&The maximal matrix $\PMmax$.\\
\bottomrule
\end{tabularx}
\end{center}
\vspace{1em}

Throughout we set $n_r$ and $n_c$ equal to, respectively, the number of rows and the number of columns of the input matrix $\PM$. A vector $s$ of length $n_c$ is used to represent the permitted permutations of the columns of $\PM$. Initially $s$ is defined as
$$s=(a,\ldots,a)\in\Z^{n_c},$$
where $a:=1+\max\PM_{ij}$ is larger than any entry of the matrix $\PM$. At each step of the recursion, the value of $n_c$ remains unchanged, but the value of $n_r$ will decrease by one as a row of $\PM$ is removed from consideration. The vector $s$ will be modified to reflect the symmetries of the previously steps; two coefficients $s_j$ and $s_k$ are equal if and only if the columns $j$ and $k$ can be exchanged without affecting the computations so far. By construction $s$ will always satisfy:
\begin{enumerate}
\item either $s_j=s_{j+1}$ or $s_j=s_{j+1} + 1$, for each $1\leq j<n_c$;
\item $s_{n_c}= a$.
\end{enumerate}

The first stage is to calculate the maximum possible row $\mmax{R}$ of $\PM$, where each row is sorted in decreasing order. Once done, we update the vector $s$ to reflect the possible column permutations that will leave $\mmax{R}$ unchanged.


\begin{algorithmic}
\State $\tilde{R}_i\gets\mathrm{Sort}_\geq\{(s_j,\PM_{ij})\mid 1\leq j\leq n_c\}$
\State $R_i\gets (\tilde{R}_{ij2}\mid 1\leq j\leq n_c)$
\State $\mmax{R}\gets \max\left\{R_i\mid 1\leq i\leq n_r\right\}$
\State $s'\gets s$
\For{$j=n_c-1$ \ForTo $1$ \By $-1$}
\If{$s_j=s_{j+1}\wedge\mmax{R_j}\ne\mmax{R_{j+1}}$}
\For{$k=1$ \ForTo $j$}
\State $s'_k\gets s'_k+1$
\EndFor
\EndIf
\EndFor
\end{algorithmic}

Next we collect together all non-isomorphic ways of writing $\PM$ with $\mmax{R}$ as the first row. These possibilities are recorded in the set $\mathcal{M}$.


\begin{algorithmic}
\State $\mathcal{M}\gets\{\}$
\For{$i=1$ \ForTo $n_r$ \SuchThat $R_i=\mmax{R}$}
\State $M\gets\mathrm{SwapRow}(\PM,1,i)$
\State $T\gets\mathrm{Sort}_\geq\{(s_j,M_{1j},j)\mid 1\leq j\leq n_c\}$
\State $\tau\gets\text{permutation in $S_{n_c}$ sending $j$ to $T_{j3}$}$
\State $M\gets\tau(M)$
\State $\tilde{M}\gets M$
\State $\tilde{M}_1\gets s'$
\If{$\bigwedge_{M'\in\mathcal{M}}\left(\tilde{M'}\not\cong\tilde{M}\right)$}
\State $\mathcal{M}\gets \mathcal{M}\cup\{M\}$
\EndIf
\EndFor
\end{algorithmic}

When all possible symmetries of the columns have been exhausted, the vector $s'$ will be equal to the sequence
$$(a+n_c-1,a+n_c-2,\ldots,a).$$
If this is the case, then $\PMmax$ is the maximum matrix in $\mathcal{M}$, once the rows have been placed in decreasing order. If there remain symmetries to explore, then we recurse on each of the matrices in $\mathcal{M}$ using the new permutation vector $s'$; $\PMmax$ is given by the largest resulting matrix.


\begin{algorithmic}
\If{$n_r=1$}
\State $\PMmax\gets\mmax{R}$
\ElsIf{$s'_1=s'_{n_c}+n_c-1$}
\State $\PMmax\gets\max\left\{\mathrm{SortRows}_\geq(M)\mid M\in\mathcal{M}\right\}$
\Else
\State $\mathcal{M'}\gets\{\}$
\For{$M\in\mathcal{M}$}
\State $M'\gets\mathrm{RemoveRow}(M,1)$
\State $M'\gets(\text{recurse with $PM\gets M'$ and $s\gets s'$})$
\State $\mathcal{M'}\gets\mathcal{M'}\cup\{M'\}$
\EndFor
\State $\PMmax\gets\mathrm{VerticalJoin}(\mmax{R},\max\mathcal{M'})$
\EndIf
\end{algorithmic}

\bibliographystyle{amsalpha}
\newcommand{\etalchar}[1]{$^{#1}$}
\providecommand{\bysame}{\leavevmode\hbox to3em{\hrulefill}\thinspace}
\providecommand{\MR}{\relax\ifhmode\unskip\space\fi MR }
\providecommand{\MRhref}[2]{%
  \href{http://www.ams.org/mathscinet-getitem?mr=#1}{#2}
}
\providecommand{\href}[2]{#2}

\end{document}